\documentclass[12pt, reqno, twoside]{amsart}
\usepackage{amsmath, amsthm, amscd, amsfonts ,, amssymb, graphicx, color}
\usepackage[bookmarksnumbered, colorlinks, plainpages]{hyperref}
\hypersetup{colorlinks=true,linkcolor=red, anchorcolor=green, citecolor=cyan, urlcolor=red, filecolor=magenta, pdftoolbar=true}
\usepackage{enumitem}
\usepackage{enumitem}
\newtheorem{theorem}{Theorem}[section]

\newtheorem{proposition}[theorem]{Proposition}
\newtheorem{corollary}[theorem]{Corollary}
\theoremstyle{definition}

\theoremstyle{remark}
\newtheorem{remark}[theorem]{Remark}
\numberwithin{equation}{section}

\begin{document}

\setcounter{page}{1}

\title[Basic equations for almost Ricci-harmonic solitons]{Basic structural  equations for almost Ricci-harmonic solitons and applications}

\author[A. Abolarinwa]{Abimbola Abolarinwa}

\address{Department of Physical Sciences
Landmark University, P. M. B. 1001, Omu-Aran, Kwara State,
Nigeria.}
\email{\textcolor[rgb]{0.00,0.00,0.84}{A.Abolarinwa1@gmail.com}}

\subjclass[2010]{53C21, 53C25, 53C44, 53C65}

\keywords{Ricci soliton, harmonic maps, Ricci harmonic soliton, rigidity}

\date{May
27, 2017}

\begin{abstract}
This paper studies gradient almost Ricci-harmonic soliton  with respect to a fixed metric. We rely on     analytic techniques to estabilish some basic elliptic and integral equations for the structure of almost Ricci-harmonic soliton  which generalizes that of Ricci-hamonic solitons on one hand and that of almost Ricci soliton  on the other hand.
The consequences of these formulas in relation to classification and rigidity  results for  gradient almost Ricci-harmonic solitons are also discussed.  
\end{abstract}
 \maketitle
 

\section{Introduction and the setting}
Let $(M, g)$ be a complete Riemannian manifold.  $(M, g)$ is said to have Ricci soliton structure if it satisfies the soliton equation, for any smooth vector field $X$ and some constant $\sigma \in \mathbb{R}$,
\begin{align}\label{eqn1}
Rc + \frac{1}{2} \mathcal{L}_Xg  = \sigma g,
\end{align}
where $Rc$ is the Ricci tensor of $M$, $g$ is the metric and $\mathcal{L}_X$ is the Lie derivative in the direction of $X$. In this case, the triple $(g, X, \sigma)$ is called Ricci soliton and said to be shrinking, steady or expanding if $\sigma >0, \sigma =0  $ or $\sigma <0$,  respectively. If $X$ is a gradient of some smooth function $f$ on $M$, that is, $X = \nabla f$, Ricci soliton equation (\ref{eqn1}) is then  written as 
\begin{align}
Rc + \text{Hess} f = \sigma g,
\end{align}
where $\text{Hess}f$ is the Hessian of function $f$. Then the triple $(g, f, \sigma)$ is called gradient Ricci soliton and  $f$ is called the potential function. If  either $X$ or $\nabla f$ is a Killing vector field  or $f$ is a constant, we say that the soliton is trivial and the underlying metric is Einstein. Compact Ricci solitons are special solutions to Hamilton's Ricci flow \cite{[Ha82]} and they occur as the flow's singularity models or blow-up limits \cite{[HDcao],[Pe02]}. We refer the reader to \cite{[HDcao],[CK04],[CLN06]} for background on Ricci solitons and their connection to the Ricci flow. Ricci soitons are  also very useful in theoretical physics as they have connection with renormalization group flow, Einstein field and string theories \cite{[CV9],[Fr85]}.  
A generalization of Ricci solitons is the Ricci-harmonic solitons which serve as special solutions of the Ricci-harmonic flow. 
\subsection{The Ricci-harmonic flow}
Let $(M, g)$ and $(N, h)$ be compact  Riemannian
manifolds of dimensions $m$ and $n$ respectively.
Let a smooth map $ u: M  \rightarrow N$ be a critical point of the  
energy integral 
 $E(u) = \int_M | \nabla u|_g^2 d \mu_g,$  
 where $N$ is 
isometrically embedded in $\mathbb{R}^d, \ d \geq n,$ by the Nash embedding
theorem. These critical points are called harmonic maps between $M$ and $N$, which are generalization of harmonic functions such as identity and constant maps. 
 The configuration $(g(x ,t), u(x, t)), t \in [0, T)$ of  a one parameter
family of Riemannian metrics $g(x, t)$ and a family of smooth maps $u(x, t)$ is defined to be Ricci-harmonic  flow if it satisfies the coupled system 
of nonlinear parabolic equations 
\begin{equation}\label{eq11}
\left. \begin{array}{l}
\displaystyle \frac{\partial}{\partial t} g(x ,t) = - 2 Rc (x ,t) + 2 \alpha(t) \nabla
u(x ,t) \otimes \nabla u(x ,t)
\\ \ \\
\displaystyle \frac{\partial}{\partial t} u(x ,t) = \tau_g u(x ,t).
\end{array} \right.
\end{equation}
  Here  $\alpha(t) \equiv 
\alpha > 0$ is  a time-dependent coupling constant  and $\tau_g u$ is the intrinsic 
Laplacian of $u$ which denotes the tension field of map $u$.  The system 
(\ref{eq11}) was first studied by  List  \cite{[Li08]} in a special case $(N,h)=(\mathbb{R}, dr^2)$,   with motivation coming from general relativity. When $u$ is a constant function, List's flow becomes the well known Hamilton's Ricci flow.
List's flow was generalized 
by  M\"uller  \cite{[Mu12]} to the general case $ N \hookrightarrow \mathbb{R}^d$,
for sufficiently large $d$ as we see in system (\ref{eq11}). Precisely, the system couples together the Hamilton's 
Ricci flow \cite{[Ha82]} and Eells and Sampson's heat flow for harmonic maps
\cite{[ES64]}. The Ricci-harmonic flow is closer to the Ricci flow in behaviours,  though   less singular.  
 For a detail introduction to the  
study of the Ricci flow see the book by Chow and Knopf \cite{[CK04]} or  Chow,
Lu and Ni \cite{[CLN06]}. The note \cite{[EL1]} will provide necessary background on harmonic map theory.  

The short-time existence, uniqueness and obstruction to long-time existence of Ricci-harmonic flow are   discussed in \cite{[Mu12]}, see also \cite{[Li16]}.  In recent time, several researchers have extended fundamental results such as Perelman's entropy formulas, no breather theorems, noncollapsing theorems,   volume growth, gradient estimates and Harnack inequalities in the study of Ricci flow to the setting of Ricci-harmonic flow, for example, see  \cite{[Ab1],[Ab2],[Ab4],[GPT2],[Li08],[MW],[Ta1],[Ta2],[Wa],[GWu]}.

\subsection{Gradient Ricci-harmonic solitons}
Let $f : M \to \mathbb{R}$ be a smooth function and $u : (M,g) \to (N,h)$ a smooth map (not necessarily harmonic map), where $(M,g)$ and $(N,h)$ are static Riemannian manifolds. Then, the quadruple $(g,u,f, \sigma)$ is called a gradient Ricci-harmonic soliton if it satisfies the coupled system of elliptic partial differential equations
\begin{equation}\label{eq111}
\left. \begin{array}{l}
\displaystyle  Rc  - \alpha \nabla
u \otimes \nabla u + \text{Hess} f  =   \sigma g 
\\ \ \\
\displaystyle  \hspace{2cm} \tau_g u  - \langle \nabla u, \nabla f\rangle = 0,
\end{array} \right.
\end{equation}
where $\alpha >0$ is a positive constant depending on $m$ and $\sigma \in \mathbb{R}$. The function $f$ is called the potential. The gradient Ricci-harmonic soliton is said to be shrinking, steady or expanding depending on whether $  \sigma >0,\  \sigma =0$ or $  \sigma<0$. We say that the soliton is trivial if $f$ is set to be constant. Here are some examples as they appear in certain contexts:
\begin{enumerate}[label=(\arabic*)] 
\item If $u$ is a constant map, then (\ref{eq111}) becomes
$$ Rc  + \text{Hess}  f   =\sigma g,$$
which is the equation of gradient Ricci soliton.
\item If $(N,h) = (\mathbb{R}, dr^2)$ and $u$, $f$ are constants, then (\ref{eq111}) reduces to
$$ Rc  = \sigma g ,$$
which says that $(M,g)$ is Einstein.  
\item If $f$ is a constant and $u$ is a harmonic map, then (\ref{eq111}) reduces to
\begin{equation*}
\left. \begin{array}{l}
\displaystyle  Rc  - \alpha \nabla
u \otimes \nabla u   =   \sigma g 
\\ \ \\
\displaystyle  \hspace{2.5cm} \tau_g u = 0,
\end{array} \right.
\end{equation*} 
which is usually called harmonic-Einstein and generalizes the notion of Einstein manifold. 
\end{enumerate}
For generalization and some extension of results on gradient Ricci solitons see \cite{[CSW11],[CZ13],[ENM],[JW13],[PRRS],[Wi]}. The paper \cite{[HDcao]} (see also the references therein) provides details on the geometry of Ricci solitons.  
Recently, a particular generalization of gradient Ricci soliton, called gradient almost Ricci soliton, was proposed in \cite{[PRRS]}. The concept appears to be natural and meaningful, though, almost Ricci soliton are not necessarily gradient. See \cite{[PRRS],[BR11],[BR12],[BR13],[FFGP]} for detail discussion on existence, rigidity, triviality, classification and so on for gradient almost Ricci soliton. 
 In this paper however,  we are considering gradient almost Ricci-harmonic solitons, in which case, $\sigma$   appearing in (\ref{eq111}) is not necessarily a constant but a smooth function on $M$. In this setting we adopt similar terminology, where we say gradient almost Ricci-harmonic soliton is shrinking, steady  or expanding  depending on whether  $\sigma$   is positive, null or negative on $M$. The gradient almost Ricci-harmonic soliton is said to be trivial if $f$  is a constant. 

 The organization of the paper is as follows: 
In Section \ref{sec2}, we give detail of notation adopted in the paper and statements of results to be proved later. Sections  \ref{sec3} discusses basic equations which will be used to prove some classification and rigidity results for almost Ricci harmonic solitons. The last section is devoted    
to the proofs of main results.
\section{Notation and Statement of results}\label{sec2} 
 \label{sec7}
\subsection{Notation}
 In the sequel, Ricci harmonic soliton  will be abbreviated as $(RH)_\alpha$ soliton. 
We write $(0,2)$-tensor $Sc := Rc - \alpha \nabla u \otimes \nabla u$,
its components as $S_{ij} := R_{ij} - \alpha \nabla_i u \nabla_j u $ in local 
coordinates while its metric trace $g^{ij}  S_{ij}$ is denoted by $S := R 
 - \alpha |\nabla u|^2$, where $Rc$, $R_{ij}$ and $R$  
 are the Ricci curvature tensor, its components
and scalar curvature, respectively. $\text{div}$ is the divergence operator, i.e., $(\text{div}   Sc)_k = g^{ij}\nabla_i S_{jk}.$
 For a smooth map $u(x) : (M, g) \rightarrow (N, h)$, where
$ g = g_{ij} dx^i dx^j$ and $h=h_{\nu \gamma } d u^\nu d u^\gamma$,  
 the energy density of $u(x)$ at a point $x \in M$ is 
$ |\nabla u|_g^2 = g^{ij} h_{ \nu \gamma} \partial_i u^\nu \partial_j u^\gamma, $
i.e., the ${\text trace}_g  u^* h$ and  $\nabla u \otimes \nabla  u = u^* h$ is the pullback of the 
metric tensor $h$ on $N$ via the map $u$. The tension field of map $u$ is denoted by 
$\tau_g u = {\text trace}_g \nabla d u$, which is a section of the bundle $u^*(T_uN)$,
where $\nabla$ is the connection on the vector bundle $T^*M \otimes u^*(TN)$ 
induced by the Levi-Civita connections on $M$ and $N$.  

For a smooth function $f$ on $M$, we write $|\nabla f|^2 = g^{ij} 
\nabla_i f \nabla_j f$ and $\Delta$ is the Laplace Beltrami operator, while the volume measure on $M$ is denoted by $ dM 
= \sqrt{|g|} dx^i, i = 1, 2, ..., m$, where $|g|$ is the determinant of $g$. 
  
\subsection{Statement of results}

 The main results of this paper which concern classification and rigidity of the gradient almost $(RH)_\alpha$ solitons  are stated in this section. Basic formulas which will help in proving them are dealt with in Section \ref{sec3}.
\begin{theorem}\label{thm21}
Let $(M,g)$ be a compact Riemannian manifold. Suppose $(g,u,f,\sigma)$ is a gradient almost $(RH)_\alpha$ soliton such that  $\Delta \sigma \leq 0$  on $M$. Setting $S_{min} := \min_M S$, $\displaystyle \sigma_* =\inf_M \sigma$ and $\displaystyle \sigma^* =\sup_M \sigma$. Then
\begin{enumerate}[label=(\roman*)] 
\item if gradient almost $(RH)_\alpha$ soliton  is steady, then  $S_{min}=0$.
\item if gradient almost $(RH)_\alpha$ soliton  is shrinking with $0<\sigma \leq \sigma^*$, then $0\leq S_{min}\leq m\sigma^*$. In particular $S\geq 0$.
\item if gradient almost $(RH)_\alpha$ soliton  is expanding with $\sigma_* \leq \sigma <0$, then $m\sigma_* \leq S_{min}\leq 0$. In particular $S\geq m\sigma$.
\end{enumerate}
\end{theorem}
\begin{corollary}\label{cor22}
With the assumptions of Theorem \ref{thm21}.
\begin{enumerate}[label=(\roman*)] 
\item In cases of shrinking and expanding (i.e., cases (ii) and (iii) of Theorem \ref{thm21}), $(M,g)$ has nonnegative scalar curvature, $R \geq \alpha |\nabla u|^2 \geq 0$.
\item In the shrinking case, if in addition $S_{min} = m \sigma^*$ then $f$ is constant, $u$ is harmonic, $S\equiv S_{min} = m \sigma^*$  and $Sc=\sigma g$. Meaning that gradient shrinking almost $(RH)_\alpha$  soliton is trivial.
\end{enumerate}
\end{corollary}
The next result is the generalization to almost $(RH)_\alpha$  soliton of \cite[Theorem 1.1]{[Wei1]} for Ricci  solitons and \cite[Theorem 3]{[BR11]} for almost Ricci solitons.
\begin{theorem}\label{thm23}
Let $(M,g)$ be compact with $m\geq 3$. Then gradient almost $(RH)_\alpha$  soliton is rigid in the sense of triviality if $$ \int_M \Big(Sc(\nabla f, \nabla f) +(m-2)\langle \nabla \sigma, \nabla f \rangle \Big) dM \leq 0,$$
\end{theorem}
where $dM$ is the volume measure of $M$.
\begin{proposition}\label{prop24}
Let $(M,g)$ be compact with $m\geq 2$. Then gradient almost $(RH)_\alpha$  soliton with $Sc \geq 0$ ($Sc \leq 0$) has constant $S$ if and only if  $Sc(\nabla f, \nabla f) = 0$ and  $\langle \nabla \sigma, \nabla f \rangle   \leq 0.$ 
\end{proposition}
If $u$ is a constant map, we obtain a version of \cite[Corollary 1]{[BR11]} as an application of integral equations proved in the next section (Theorem \ref{thm26}).
\begin{proposition}\label{prop25}
Let $(M,g)$ be compact and $(g,u,f,\sigma)$ be a nontrivial gradient almost $(RH)_\alpha$  soliton. Suppose $m \geq 2$ and $u$ is a constant map, then $(M,g)$ is isometric to an $m$-Euclidean sphere, if and only if any of the following conditions holds:
\begin{enumerate}[label=(\roman*)]  
\item $(M,g)$ has constant scalar curvature.
\item $\displaystyle  \int_M  \langle   Sc, \Delta f \rangle   dM = 0.$
\item $\displaystyle  \int_M \Big(Sc(\nabla f, \nabla f) +(m-2)\langle \nabla \sigma, \nabla f \rangle \Big) dM \leq 0.$
\item $\nabla f$ is a conformal vector field.
\end{enumerate}
Here  $dM$ is the volume measure of $M$.
\end{proposition}
We notice also that if $u$ is a nonconstant map, one can still deduce the last result. See the discussion in Remark \ref{rmk44}.

 \section{Basic equations for almost Ricci-harmonic solitons}\label{sec3} 
 Firstly, we shall prove some basic equations which will be applied to derive some classification and rigidity results. The equations are extension to almost $(RH)_\alpha$ solitons of the classical ones derived for almost Ricci solitons. We include the proofs for completeness sake and since some new terms come in.
\begin{proposition}\label{prop31}
Let $(g,u,f,\sigma)$ be a gradient almost $(RH)_\alpha$ soliton, then we have the following equations:
\begin{align}
\text{div}  Sc & = \frac{1}{2} \nabla S - \alpha \tau_g(u)\nabla u\\
\langle \nabla S, \nabla f \rangle & = 2(m-1)\langle \nabla \sigma, \nabla f \rangle + 2 Sc(\nabla f, \nabla f)\\
Sc(\nabla f, \cdot) &= \frac{1}{2} \nabla S - (m-1)\nabla \sigma\\
\nabla (S+|\nabla f|^2)  &= 2(m-1)\nabla \sigma + 2 \sigma \nabla f\\ 
\frac{1}{2} \Delta |\nabla f|^2 &= |\text{Hess}  f|^2 -(m-2)\langle \nabla \sigma, \nabla f\rangle -  Sc(\nabla f, \nabla f) + \alpha |\langle \nabla u, \nabla f|^2.
\end{align} 
\end{proposition}
The proofs of the above equations follow from standard computation which is usually performed  with the method of the moving frame in a local orthonormal frame.
\proof
The first formula is a consequence of contracted second Bianchi identity $\text{div}  Rc = \frac{1}{2} \nabla R$ and the identity $ \text{div} (\nabla u \otimes \nabla u) = \tau_g(u)\nabla u + \frac{1}{2}\nabla|\nabla u|^2$ proved in \cite[Appendix]{[GPT2]}
\begin{align*}
\frac{1}{2} \nabla R &= \text{div}  Rc\\
 & = \text{div} Sc + \alpha div(\nabla u \otimes \nabla u)\\
 & = \text{div}  Sc + \alpha \tau_g(u)\nabla u + \frac{\alpha}{2}\nabla|\nabla u|^2, 
\end{align*}
which implies
$ \frac{1}{2} \nabla(R+\alpha\nabla|\nabla u|^2) = \text{div}  Sc + \alpha \tau_g(u)\nabla u $, the first formula to be proved.

Note that the first equation in (\ref{eq111}) implies
\begin{align}\label{el1}
S_{ij}=\sigma g_{ij} - f_{ij}.
\end{align}
Taking covariant derivative of (\ref{el1}) we have 
\begin{align}\label{el2}
S_{ij,k}=\sigma_k g_{ij} - f_{ijk} 
\end{align}
since $g_{ij,k} =0$ (covariant derivative of $g$) by metric compatibility of Levi-Civita connection. Tracing with respect to $j$ and $k$ we have
\begin{align}\label{el3}
S_{ik,k}=\sigma_i   - f_{ikk}.
\end{align}
Recall the Ricci identity
\begin{align}\label{el4}
f_{ikk} - f_{kki} =R_{ik}f_k
\end{align}
and the first formula in the proposition
\begin{align}\label{el5}
S_{ij,i}= \frac{1}{2} S_j - \alpha \tau_g(u)\nabla u.
\end{align}
Now, putting (\ref{el4}) and (\ref{el5}) into (\ref{el3}), we have 
\begin{align}\label{el6}
\frac{1}{2}S_i - \alpha \tau_g(u)\nabla u = \sigma_i -f_{kki} -R_{ik}f_k. 
\end{align}
Tracing (\ref{el2}) with respect to $i$ and $j$, we have 
\begin{align}\label{el7}
 S_i   = m \sigma_i -f_{kki}.  
\end{align}
Substituting (\ref{el7}) into (\ref{el6}) we have
\begin{align}\label{el8}
 S_i   = 2(m-1) \sigma_i + 2R_{ik}f_k - 2\alpha \tau_g(u)\nabla u.  
\end{align}
In particular, we have from (\ref{el8})
\begin{align}\label{el9}
 \left. \begin{array}{l}  
\displaystyle \langle\nabla S, \nabla f\rangle = 2(m-1) \langle\nabla \sigma, \nabla f\rangle + 2Rc(\nabla f, \nabla f) - 2 \alpha \tau_g(u)  \langle\nabla u, \nabla f\rangle \\ \ \\ 
\displaystyle \hspace{1.8cm} = 2(m-1) \langle\nabla \sigma, \nabla f\rangle + 2Sc(\nabla f, \nabla f),  
\end{array} \right.
\end{align}
where we have used $\tau_g(u) =  \langle\nabla u, \nabla f\rangle$ and $Rc = Sc + \alpha \nabla u \otimes \nabla u$. This proves the second formula in the proposition.

The third equation in the proposition can be seen as a direct consequence of (\ref{el8}) or (\ref{el9}). But we prove it this way: Taking divergence of the first equation in (\ref{eq111}) yields
\begin{align}\label{el10}
{\text div }   Rc - \alpha {\text div}(\nabla u \otimes \nabla u) +    \text{div Hess} f = \text{div} (\sigma g).
\end{align}
Routine computation gives $ \text{div} (\nabla u \otimes \nabla u) = \tau_g(u)\nabla u + \frac{1}{2}\nabla|\nabla u|^2$ and
 $ {\text div Hess} f = \nabla \Delta f + Rc(\nabla f, \cdot) = m \nabla \sigma - \nabla S + Rc(\nabla f, \cdot)$. Then (\ref{el10}) implies 
\begin{align}\label{el11}
{\text div }   Rc  +(m-1)\nabla \sigma + Rc(\nabla f, \cdot) - \alpha \tau_g(u) \nabla u - \nabla S - \frac{\alpha}{2}\nabla|\nabla u|^2=0.
\end{align}
Using ${\text div} Rc = \frac{1}{2} \nabla R$ in (\ref{el11}) we have 
\begin{align}\label{el11b}
\frac{1}{2} \nabla (R - \alpha |\nabla u|^2) - \nabla S   +(m-1)\nabla \sigma + (Rc - \alpha \nabla u \otimes \nabla u)(\nabla f, \cdot) =0.
\end{align}
Then
\begin{align} 
-\frac{1}{2} \nabla S + (m-1)\nabla \sigma + Sc(\nabla f, \cdot) =0,
\end{align}
which is the third equation of the proposition.

The fourth equation is another consequence of (\ref{el8}). It has been previously derived in \cite{[BR11],[FFGP]} by different methods. Equation (\ref{el8}) implies
\begin{align*}
\nabla S & = 2(m-1) \nabla \sigma + 2 Sc(\nabla f, \cdot) \\
& =  2(m-1) \nabla \sigma + 2(\sigma g - \text{Hess} f) \nabla f \\
&=  2(m-1) \nabla \sigma + 2 \sigma \nabla f - 2 \text{Hess} f(\nabla f)\\
&= 2(m-1) \nabla \sigma + 2 \sigma \nabla f - \nabla |\nabla f|^2,
\end{align*}
which implies
\begin{align}
\nabla (S + |\nabla f|^2) - 2(m-1) \nabla \sigma - 2\sigma \nabla f = 0.
\end{align}
The fourth  equation is proved.

Inserting $\Delta f = m \sigma -S$ into the well known  Bochner formula we have
\begin{align*}
\frac{1}{2} \Delta |\nabla f|^2 & = |\text{Hess} f|^2 + \langle \nabla \Delta f, \nabla f \rangle + Rc(\nabla f, \nabla f)\\
& = |\text{Hess} f|^2 + m \langle \nabla \sigma, \nabla f \rangle - \langle \nabla S, \nabla f \rangle + Rc(\nabla f, \nabla f).
\end{align*}
Using the second equation of the proposition yields
\begin{align*}
\frac{1}{2} \Delta |\nabla f|^2 & = |\text{Hess} f|^2 + m \langle \nabla \sigma, \nabla f \rangle -2 (m-1) \langle \nabla \sigma, \nabla f \rangle - 2 Sc(\nabla f, \nabla f) \\
& \hspace{1cm} + Rc(\nabla f, \nabla f)\\
& = |\text{Hess} f|^2 -(m-2) \langle \nabla \sigma, \nabla f \rangle  - Sc(\nabla f, \nabla f) + \alpha \nabla u \otimes \nabla u (\nabla f, \nabla f),
\end{align*}
which implies the fifth equation.
Therefore the proof of the proposition is complete.

\qed

Lastly, we derive some integral equations  which will be used to derive an application which generalizes results obtained in \cite{[BR11]} for almost Ricci solitons.
\begin{theorem}\label{thm26}
Let $(M,g)$ be   compact  and  $(g,u,f,\sigma)$ be gradient almost $(RH)_\alpha$ soliton. Then we have
\begin{align}\label{eq37b} 
\displaystyle \int_M \Big| \text{Hess} f  - \frac{1}{m} g \Delta f\Big|^2 dM = \frac{m-2}{2m}\int_M \langle \nabla S, \nabla f\rangle dM - \alpha \int_M |\langle \nabla u, \nabla f \rangle|^2dM  
\end{align}
and
\begin{align}\label{eq37c}
 \left. \begin{array}{l}  
\displaystyle \int_M \Big|\text{Hess} f  - \frac{1}{m} g \Delta f\Big|^2 dM = \frac{m-2}{m}\int_M \Big(Sc(\nabla f, \nabla f) + (m-1)\langle \nabla \sigma, \nabla f\Big) dM \\ 
\displaystyle \hspace{5cm} - \alpha \int_M |\langle \nabla u, \nabla f\rangle|^2 dM,  
\end{array} \right.
\end{align}
where $dM$ is the volume measure of $M$.
\end{theorem} 
\proof
Taking divergence of the fourth equation in Proposition \ref{prop31} yields
\begin{align}\label{e1}
\Delta S + \Delta |\nabla f|^2 -2(m-1)\Delta \sigma - 2 \langle \nabla \sigma, \nabla f \rangle - 2\sigma \Delta f = 0.
\end{align}
Using second and fifth equations in Proposition \ref{prop31} into (\ref{e1}) we have 
\begin{align}\label{e2}
\Delta S + 2 |\text{Hess} f|^2 - \langle \nabla S, \nabla f \rangle  -2(m-1)\Delta \sigma   - 2\sigma \Delta f + 2\alpha|\langle \nabla u, \nabla f\rangle|^2 = 0.
\end{align}
By using the identity
$$\Big|\text{Hess} f - \frac{1}{m}g \Delta f \Big|^2 =|Hess f|^2 - \frac{1}{m} (\Delta f)^2$$
(\ref{e2}) can hen be written as 
\begin{align*} 
\Delta S +  2\Big|\text{Hess} f - \frac{1}{m}g \Delta f \Big|^2 - \langle \nabla S, \nabla f \rangle  -2(m-1)\Delta \sigma \\
   - 2 \Delta f \Big(\sigma - \frac{1}{m}\Delta f\Big) + 2\alpha|\langle \nabla u, \nabla f\rangle|^2 = 0.
\end{align*}
Note that $\frac{1}{m}S =  \sigma - \frac{1}{m}\Delta f$ (from $\Delta f = m\sigma -S$), then 
\begin{align*} 
\Delta S +  2\Big|\text{Hess} f - \frac{1}{m}g \Delta f \Big|^2 - \langle \nabla S, \nabla f \rangle  -2(m-1)\Delta \sigma \\
  - \frac{2}{m}  S\Delta f   + 2\alpha|\langle \nabla u, \nabla f\rangle|^2 = 0.
\end{align*}
Integrating over $M$ and using compactness of $M$, we have
\begin{align}\label{e3} 
\displaystyle   2 \int _M \Big|\text{Hess} f - \frac{1}{m}g \Delta f \Big|^2 dM  = \frac{m-2}{m} \langle \nabla S, \nabla f \rangle  
 - 2\alpha \int _M|\langle \nabla u, \nabla f\rangle|^2  dM
\end{align}
from where (\ref{eq37b}) follows.  Now using the second equation of Proposition \ref{prop31} in (\ref{e3}) we obtain (\ref{eq37c}) at once. The proof is complete.

\qed

 \section{Proofs of results}
The aim of this section is to use some of the equations obtained in the last  section to classify gradient almost $(RH)_\alpha$ solitons and determine when they are rigid as stated in Section \ref{sec2}. The first result to prove is Theorem \ref{thm21}.

\subsection{Proof of Theorem \ref{thm21}}
To start with, taking divergence of the fourth equation in Proposition \ref{prop31} gives
\begin{align}\label{e41}
\Delta S + \Delta |\nabla f|^2 -2(m-1)\Delta \sigma - 2 \langle \nabla \sigma, \nabla f \rangle - 2\sigma \Delta f = 0.
\end{align}
Using  fifth  and second equations in Proposition \ref{prop31}, identity  $\Delta f = m\sigma -S$, the condition $\Delta \sigma \leq 0$  and the inequality $|Sc|^2 \geq \frac{2}{m} S^2$ into (\ref{e41}) we have 
\begin{align*} 
0 &= \Delta S + 2 |\text{Hess} f|^2 - \langle \nabla S, \nabla f \rangle  -2(m-1)\Delta \sigma   - 2\sigma \Delta f + 2\alpha|\langle \nabla u, \nabla f\rangle|^2 \\
& = \Delta S + 2 |\sigma g - Sc|^2 - \langle \nabla S, \nabla f \rangle  -2(m-1)\Delta \sigma   - 2\sigma (m \sigma-S) + 2\alpha|\langle \nabla u, \nabla f\rangle|^2  \\
& = \Delta S + 2 |Sc|^2 -2 \sigma S - \langle \nabla S, \nabla f \rangle  -2(m-1)\Delta \sigma + 2\alpha|\langle \nabla u, \nabla f\rangle|^2 \\
&\geq \Delta S + \frac{2}{m}S^2 - 2 \sigma S - \nabla S(\nabla f).
\end{align*}
Since $M$ is compact, there exists a point $q \in M$, where $S$ attains its minimum. By the maximum principle, at this point $\nabla S(q) = 0$ and $\Delta S(q) \geq 0$. Then the last inequality implies
\begin{align}\label{e42}
S_{min}(S_{min} - m \sigma) \leq 0.
\end{align} 
The first case is when $\sigma =0$, that is, the case of steady gradient almost $(RH)_\alpha$ soliton, we have $S_{min} =0$ by (\ref{e42}).
In the case of shrinking gradient almost $(RH)_\alpha$ soliton $\sigma >0$ with $0<\sigma \leq \sigma^*$, by (\ref{e42}) we have $0\leq S_{min} \leq m \sigma^*$.
In the case of expanding gradient almost $(RH)_\alpha$ soliton $\sigma <0$ with $\sigma_* \leq \sigma<0$,  we have $m \sigma_*  \leq S_{min}<0$ and $S \geq S_{min} \geq m \sigma_*$.

\qed
 
 \begin{remark}
 If one follows the steps in the proof of Theorem 0.4 of \cite{[PRRS]}, it can be further shown that under a pairwise control on $\sigma$,  the curvature function $S$ of gradient almost $(RH)_\alpha$ soliton is bounded from below and the lower bound can be estimated from both below and above.
 \end{remark}
We note that when $S$ attains its minimum, by the maximum principle   the inequality
\begin{align}\label{eqs}
\Delta S_{min} + \frac{2}{m}S_{min}^2 - 2 \sigma S_{min} - \nabla S_{min}(\nabla f) \leq 0
\end{align}
implies that if $S$ is nonconstant, then it must be everywhere positive and then $\sigma$ is positive, which says that gradient almost $(RH)_\alpha$ soliton is shrinking. This implies the following proposition.
\begin{proposition}
Every steady or expanding gradient almost $(RH)_\alpha$ soliton has constant $S$.
\end{proposition}
\proof
By the weak maximum principle inequality (\ref{eqs}) implies
\begin{align} 
\frac{1}{2}\Delta S_{min} \leq \frac{1}{m}S_{min}(m \sigma_* - S_{min}).  
\end{align}
This implie that $\Delta S_{min} \leq 0$ since $S_{min}=0$ ( for the steady case) and $S_{min} \geq m \sigma_*$ (for the expanding case). Then $ S \geq S_{min} \geq 0$ is nonnegative superharmonic and it is therefore constant.

\qed

The above proposition forces $\Delta f = m \sigma - S$ to be $\Delta f \geq 0$ since $S = m \sigma_* \leq m \sigma$. Thus $f$ is subharmonic and the compactness of $M$ implies that $f$ is constant and the soliton is trivial. Hence, the following corollary.
\begin{corollary}
Every steady or expanding gradient almost $(RH)_\alpha$ soliton  ($S= m \sigma_*$) is trivial on compact $M$.
\end{corollary}

\subsection{Proof of Theorem \ref{thm23}}
Integrating both sides of the fifth equation in Proposition \ref{prop31}, we have 
\begin{align*}
\frac{1}{2} \int_M \Delta |\nabla f|^2 dM &= \int_M |\text{Hess} f|^2 dM - \int_M Sc(\nabla f, \nabla f) dM \\
& - (m-2) \int_M \langle \nabla \sigma, \nabla f\rangle dM + \alpha \int_M |\langle \nabla u, \nabla f\rangle|^2 dM,
\end{align*}
where $dM$ is the volume measure of $M$. Since $M$ is compact we have 
\begin{align}\label{e43}
 \int_M |\text{Hess} f|^2 dM  & =  \int_M \Big(Sc(\nabla f, \nabla f)  + (m-2) \langle \nabla \sigma, \nabla f\rangle\Big) dM \\
 &-  \alpha \int_M |\langle \nabla u, \nabla f\rangle|^2 dM.
\end{align}
Note that the second on the right hand side of (\ref{e43}) is nonnegative and by assumption the first term is nonpositive. Then we obtain $Hess f =0$ which implies $\nabla f$ is a killing vector field and that $f$ is a constant by compactness of $M$. 
\qed

\subsection{Proof of Proposition \ref{prop24}}
The proof of Proposition \ref{prop24} can be established using any of the first three equations in Proposition \ref{prop31}. In what follows, we mimick the proof of Propositon 3 in \cite{[Wei1]}.
\proof
Consider a nonnegative (nonpositive) definite, self-adjoint operator $T$, we know that
 $\langle Tv, v \rangle =0$ implies $Tv=0$. Taking $T$ to be $(0,2)$-tensor, the proposition follows from the second or third equation in Proposition \ref{prop31}. That is,
$$\frac{1}{2} \nabla S (\nabla f) =  (m-1)\langle \nabla \sigma, \nabla f \rangle +   Sc(\nabla f, \nabla f)$$
which implies $\nabla S = 0$ by the assumption. 

\subsection{Proof of Propostion \ref{prop25}}
\proof
Considering any of the conditions in the proposition, we have from either (\ref{eq37b}) or (\ref{eq37c}) that 
\begin{align}\label{eq46}
\int_M \Big| \text{Hess} f  - \frac{1}{m} g \Delta f\Big|^2 dM = 0. 
\end{align}
Indeed, supposing condition (ii), we have 
\begin{align*}
\displaystyle 0 = \int_M \langle Sc, \Delta f \rangle dM & = \int_M  \langle \text{div} Sc, \nabla f \rangle dM \\
\displaystyle & = - \frac{1}{2} \int_M  \langle \nabla S, \nabla f \rangle dM + \alpha \int_M |\langle \nabla u, \nabla f \rangle| dM,
\end{align*}
which holds if and only if $S$ is scalar, since $u$ is contant and $\nabla f$ is not killing (nontrivial soliton). Note that (\ref{eq46}) implies $Sc - \frac{1}{m}Sg =0$ (see \ref{eq47c} below) and consequently, $Rc = \frac{1}{m}Rg$, thus, $(M,g)$ is Einstein. Therefore (\ref{eq11}) yields $Hess f = (\sigma - \frac{1}{m} R)g$, saying that $\nabla f$ is a nontrivial conformal vector field.

If on the other hand, we suppose $\nabla f$ is a conformal vector field, then (\ref{eq37b}) will yield $\int_M |\langle \nabla S, \nabla f \rangle| dM \geq 0$. The fact that $u$ is constant and (\ref{eq46}) give equality $\int_M |\langle \nabla S, \nabla f \rangle| dM = 0$, meaning that $(M, g)$ has constant curvature.  

The rest of the proof is now similar to that of \cite[Corollary  1]{[BR11]}, where a theorem of Tashiro \cite{[Tash]} together with Theorem \cite[Theorem 2]{[BR11]}  was applied. The rest argument is therefore omitted here. Tashiro's theorem says that a complete Riemannian manifold $M^m$, $m\geq 2$ is conformally equivalent to the Euclidean sphere of the same dimension if it admits a nontrivial solution of the equation $\text{Hess} \rho - \frac{1}{m}g \Delta \rho = 0$, where $\rho$ is a smooth function on $M$. \cite[Theorem 2]{[BR11]} says that when $m\geq 3$,  $(M, g)$ is  isometric to Euclidean sphere $\mathbb{S}^m$ if $u$ is a constant and $\nabla f$ is a nontrivial conformal field.

\qed

\begin{remark}\label{rmk44}
If $u$ is a nonconstant map, one can still deduce the conclusion of Propostion \ref{prop25}. Noticing that (\ref{eq37b}) and (\ref{eq37c}) under Theorem \ref{thm26} respectively become (since $\alpha >0$) 
\begin{align}\label{eq45c} 
\displaystyle \int_M \Big| \text{Hess} f  - \frac{1}{m} g \Delta f\Big|^2 dM \leq \frac{m-2}{2m}\int_M \langle \nabla S, \nabla f\rangle dM   
\end{align}
and
\begin{align}\label{eq46c}  
\displaystyle \int_M \Big| \text{Hess} f  - \frac{1}{m} g \Delta f\Big|^2 dM \leq \frac{m-2}{m}\int_M \Big(Sc(\nabla f, \nabla f) + (m-1)\langle \nabla \sigma, \nabla f\Big) dM. 
\end{align}
\end{remark}
One can then invoke any of the conditions (i)--(iv) of Propostion \ref{prop25} and conlude that 
\begin{align*}  
\displaystyle \int_M \Big| Sc  - \frac{1}{m} Sg\Big|^2 dM \leq   0  
\end{align*}
since
\begin{align}\label{eq47c} 
\text{Hess} f  - \frac{1}{m} g \Delta f = \text{Hess} f - (\sigma - \frac{1}{m}S)g = -Sc + \frac{1}{m} Sg.
\end{align} 
 Therefore    
\begin{align*} 
0 & = Sc  - \frac{1}{m} Sg   = Rc - \frac{1}{m} (R + \alpha (m-1)  |\nabla u|^2) g 
 \end{align*}
and we conclude that $(M,g)$ is  Einstein, since every Einstein manifold belongs to the
class of Riemannian manifolds $(M, g)$ realizing the relation $R_{ij} = a g_{ij} + b \nabla_i v \nabla_j v $, 
where $a$ and $b$ are certain nonzero constants (see \cite{[CM20]}).
 
  \section*{Acknowledgement}
 The author wishes to thank the anonymous referees for their useful comments. 
 

\end{document}